\documentclass[12pt, twoside, leqno]{article}

\usepackage{amsmath,amsthm}
\usepackage{amssymb}

\usepackage{enumitem}

\usepackage[T1]{fontenc}
\usepackage[latin9]{inputenc}
\usepackage{geometry}
\usepackage{mathrsfs}
\usepackage{amsmath}
\usepackage{amssymb}
\usepackage{stackrel}
\usepackage{setspace}
\PassOptionsToPackage{normalem}{ulem}
\usepackage{ulem}
\makeatletter
\@ifundefined{date}{}{\date{}}
\makeatother

\usepackage{babel}
\begin{document}

\baselineskip=17pt

\title{The Return Times Theorem, Auto-Correlation and Sequences with an
Empty Fourier-Bohr Spectrum}

\author{Matan Tal\\
The Hebrew University of Jerusalem}
\maketitle

\renewcommand{\thefootnote}{}



\renewcommand{\thefootnote}{\arabic{footnote}}
\setcounter{footnote}{0}

\begin{abstract}
This paper explores the proof by J. Bourgain, H. Furstenberg, Y. Katznelson
and D.S. Ornstein of their return times theorem \cite{key-3} and
lights a corner in it regarding the role of auto-correlation. As for
pointwise convergence, this was already observed in \cite{key-9},
and here we exploit the opportunity to write down the proof. This
yields a more intrinsic characterization of the sequences satisfying
the pointwise theorem. Then we proceed and obtain a characterization
linked to auto-correlation also to sequences satisfying the mean theorem
- by that theorem those were already known to be exactly the sequences
with an empty Fourier-Bohr spectrum. Some further investigation is
done and examples are provided regarding generic sequences satisfying
the pointwise theorem for which the measure on the circle that the
auto-correlation function represents (by Fourier transform) is not
atomless, and also regarding the existence of sequences that satisfy
the mean theorem but not the pointwise one.\\
\end{abstract}

\section{Introduction}

Throughout, we work in $\mathbb{C}^{\mathbb{N}}$. We denote its shift
map by $\sigma$, and the function that gives the value at the first
index by $\pi_{1}:\mathbb{C}^{\mathbb{N}}\rightarrow\mathbb{C}$.\\

In this paper, our matter of interest is sequences $x=\left(x_{n}\right)_{n\in\mathbb{N}}\in\mathbb{C}^{\mathbb{N}}$
that satisfy $\frac{1}{T}\stackrel[n=1]{T}{\sum}x_{n}Y_{n}\left(\omega\right)\xrightarrow[N\rightarrow\infty]{}0$
almost-surely for any complex-valued stationary process $\left(Y_{n}\left(\omega\right)\right)_{n\in\mathbb{Z}}$
of finite variance - we shall call them \textit{pointwise cancellation sequences} - and also sequences for
which  \\  $\left\Vert \frac{1}{T}\stackrel[n=1]{T}{\sum}x_{n}Y_{n}\left(\omega\right)\right\Vert _{2}\xrightarrow[T\rightarrow\infty]{}0$
for every such process that we shall call \textit{mean cancellation
sequences}. \\

By a short argument that will appear later in the paper (more specifically,
in the proof of Theorem 2) we shall see that $x$ is a mean cancellation
sequence if and only if $\left\Vert \frac{1}{T}\stackrel[n=1]{T}{\sum}x_{n}Y_{n}\left(\omega\right)\right\Vert _{2}\xrightarrow[T\rightarrow\infty]{}0$
for every \uline{bounded} complex-valued stationary process $\left(Y_{n}\left(\omega\right)\right)$.
So for $x$ satisfying \\ $\sup_{N\in \mathbb{N}} \frac{\left|x_{1}\right|+\dots+\left|x_{N}\right|}{N}<\infty$
(we shall anyway be concerned only with such $x$), pointwise cancellation
implies mean cancellation by Lebesgue dominated convergence theorem.
In Section 5 we shall see that there exist bounded sequences that
are mean cancellation sequences but not pointwise ones.\\

In the opening paragraph of \cite{key-3} the following result is
proved with a direct application of the spectral theorem on the
unitary shift map ($Y_{n}\mapsto Y_{n+1}$) and the vector $Y_{1}$.\\

\textbf{Theorem 0:} \textit{Let $x=\left(x_{n}\right)_{n\in\mathbb{N}}$
be a complex sequence satisfying that $\sup_{N\in \mathbb{N}} \frac{\left|x_{1}\right|+\dots+\left|x_{N}\right|}{N}<\infty$.
The following conditions are equivalent:}

\textbf{\textit{(i) }}\textit{$\frac{1}{T}\stackrel[n=1]{T-1}{\sum}x_{n}z^{n}\xrightarrow[T\rightarrow\infty]{}0$
for all $z\in\mathbb{C}$ on the unit circle.}

\textbf{\textit{(ii)}}\textit{ $x=\left(x_{n}\right)_{n\in\mathbb{N}}$
is a mean cancellation sequence.}\\

(The implication (ii)$\rightarrow$(i) is not mentioned in \cite{key-3}
but its proof is an immediate conclusion from the fact that for every
$z_{0}$ on the unit circle there exists stationary processes $\left(Y_{n}\left(\omega\right)\right)_{n\in\mathbb{Z}}$
for which $z_{0}$ is in the point spectrum of $Y_{1}$ relative to
the unitary shift map - for example, ones originating from circle rotations.)\\

A sequence satisfying (i) is customarily called a sequence with an
\textit{empty Fourier-Bohr spectrum}. This is a special case of sequences
$x$ for which $\lim_{T\rightarrow\infty}\frac{1}{T}\stackrel[n=1]{T-1}{\sum}x_{n}z^{n}$
exists for all $z\in\mathbb{C}$ on the unit circle, which are customarily
called \textit{Hartman almost-periodic sequences} (for instance cf.
\cite{key-5}), and in \cite{key-6} J.P. Kahane proved that for this
latter type of sequences there is at most a countable set of $z$
for which the limit is not zero. Theorem 0 gave another characterization
of sequences that have an empty Fourier-Bohr spectrum.\\

The main result in \cite{key-3} is the following theorem.\\

\textbf{Theorem 1':} \textit{Let $x=\left(x_{n}\right)_{n\in\mathbb{N}}$
be a complex sequence generic relative to a $\sigma$-invariant Borel
probability measure $\mu$ on $\mathbb{C}^{\mathbb{N}}$ satisfying \\\\
$\sup_{N\in \mathbb{N}} \frac{\left|x_{1}\right|^2+\dots+\left|x_{N}\right|^2}{N}<\infty$.
The following conditions are equivalent:}

\textbf{\textit{(i) }}\textit{$\frac{1}{T}\stackrel[n=1]{T-1}{\sum}x_{n}\overline{\xi_{n}}\xrightarrow[T\rightarrow\infty]{}0$
for $\xi\in\mathbb{C}^{\mathbb{N}}$ $\mu$-almost-surely.}

\textbf{\textit{(ii)}}\textit{ $x=\left(x_{n}\right)_{n\in\mathbb{N}}$
is a pointwise cancellation sequence.}\\

\textbf{Remark: }It was actually formulated in \cite{key-3} for a
bounded $x$ and the cancellation property refered there to processes of just finite expectation.
The formulation given above fits better the character of our paper
and requires only a straight-forward modification of the original
proof. In fact, also in \cite{key-3} the proof that (i)$\rightarrow$(ii)
is written explicitly only for processes $\left(Y_{n}\left(\omega\right)\right)_{n\in\mathbb{Z}}$
of finite variance. The deduction from that to finite expectation
requires a step which is implicit in that proof, and its details can
be found in Lemma 5 and Lemma 8 in \cite{key-2}. \\

Unlike in Theorem 0, condition (i) of Theorem 1' refers to the $\mu$
for which $x$ is generic. It
is observed in \cite{key-9} (cf. also p. 15 in \cite{key-7}) that
the proof of Theorem 1' can actually yield an equivalent condition to (ii) that, unlike (i), does not refer to $\mu$
and hence is a more intrinsic characterization of such an $x$. It
involves an expression linked to the auto-correlation of the sequence
$x$. This is Theorem $1$ in the present paper. \\

\textbf{Theorem 1: }\textit{Let $x=\left(x_{n}\right)_{n\in\mathbb{N}}$
be a complex sequence generic relative to a $\sigma$-invariant Borel
probability measure $\mu$ on $\mathbb{C}^{\mathbb{N}}$ satisfying \\\\
$\sup_{N\in \mathbb{N}} \frac{\left|x_{1}\right|^2+\dots+\left|x_{N}\right|^2}{N}<\infty$.
The following conditions are equivalent:}

\textbf{\textit{(i) }}\textit{$\frac{1}{T}\stackrel[n=1]{T-1}{\sum}x_{n}\overline{\xi_{n}}\xrightarrow[T\rightarrow\infty]{}0$
for $\mu$-almost-every $\xi\in\mathbb{C}^{\mathbb{N}}$.}

\textbf{\textit{(ii)}}\textit{ For any $\varepsilon>0$ there exists
$N'>0$ such that for any $N''>N'$ there exists $T'>0$ for which }

\textit{$\frac{\left|\left\{ 0<\tau\leq T\,:\,\exists N'\leq N\le N''\,\left|\frac{1}{N}\stackrel[n=1]{N}{\sum}x_{n+\tau}\overline{x_{n}}\right|\geq\varepsilon\right\} \right|}{T}<\varepsilon$
for every $T>T'$.}

\textbf{\textit{(iii)}}\textit{ $x=\left(x_{n}\right)_{n\in\mathbb{N}}$
is a pointwise cancellation sequence.}\\

\textbf{Remark:} The proof of the implication (ii)$\rightarrow$(iii)
(Prop. 2 below) does not make use of the fact that $x$ is generic
relative to some measure. Therefore, (ii) is a sufficient condition
for (iii) also without that assumption. We do not know whether, in
this more general situation, it is also a necessary condition or not.\\

In the next section, we write down the full proof of Theorem 1 (since
in \cite{key-9} the theorem was merely formulated).\\

Based on a proof similar in its core to the proof of Theorem 1, we
continue and prove the following theorem reinforcing Theorem 0 when \\\\
$\sup_{N\in \mathbb{N}} \frac{\left|x_{1}\right|^2+\dots+\left|x_{N}\right|^2}{N}<\infty$.\\

\textbf{Theorem 2:} \textit{Let $x=\left(x_{n}\right)_{n\in\mathbb{N}}$
be a complex sequence satisfying \\ $\sup_{N\in \mathbb{N}} \frac{\left|x_{1}\right|^2+\dots+\left|x_{N}\right|^2}{N}<\infty$.
The following conditions are equivalent:}

\textbf{\textit{(i) }}\textit{$\frac{1}{T}\stackrel[n=1]{T-1}{\sum}x_{n}z^{n}\xrightarrow[T\rightarrow\infty]{}0$
for all $z\in\mathbb{C}$ on the unit circle.}

\textbf{\textit{(ii)}}\textit{ For any increasing sequence of integers
$T_{k}$ there exists a sub-sequence $T_{k_{l}}$ that satisfies that
for any $\varepsilon>0$ there exists $r'>0$ such that for any $r''>r'$
there exists $l'>0$ for which}

\textit{$\frac{\left|\left\{ 0<\tau\leq T_{k_{l}}\,:\,\exists r'\leq r\le r''\,\left|\frac{1}{T_{k_{r}}}\stackrel[n=1]{T_{k_{r}}}{\sum}x_{n+\tau}\overline{x_{n}}\right|\geq\varepsilon\right\} \right|}{T_{k_{l}}}<\varepsilon$
for every $l>l'$.}

\textbf{\textit{(iii)}}\textit{ $\left\Vert \frac{1}{T}\stackrel[n=1]{T}{\sum}x_{n}Y_{n}\left(\omega\right)\right\Vert _{2}\xrightarrow[T\rightarrow\infty]{}0$
for any complex-valued }\textit{\uline{bounded}}\textit{ stationary
process $\left(Y_{n}\left(\omega\right)\right)_{n\in\mathbb{Z}}$.}

\textbf{\textit{(iv)}}\textit{ $x=\left(x_{n}\right)_{n\in\mathbb{N}}$
is a mean cancellation sequence.}\\

So Theorem 2 gives another characterization for a sequence to have
an empty Fourier-Bohr spectrum, this time applying an expression linked
to its auto-correlation. It is proved in section 3. To make this paper
self-contained we include in the proof also the part that has already
appeared in \cite{key-3}.\\

Turning to another line of inquiry, it was proved in \cite{key-3}
that for all ergodic stationary processes $\left(Y_{n}\left(\omega\right)\right)_{n\in\mathbb{Z}}$
of finite variance with $Y_{1}$ having a continuous spectrum relative
to the unitary shift map,
the event of being generic and also a pointwise cancellation sequence
is of probability $1$. A generic $x$ represents a stationary process
with such a continuous spectrum if and only if \[ \lim_{T\rightarrow\infty}\frac{1}{T}\stackrel[\tau=0]{T-1}{\sum}\lim_{N\rightarrow\infty}\left|\frac{1}{N}\stackrel[n=1]{N}{\sum}x_{n+\tau}\overline{x_{n}}\right|=0. \]
Section 4 is devoted to the possibility of a pointwise cancellation
sequence that is generic relative to a stationary process for which
the above mentioned spectrum is not continuous. A posteriori, after
establishing Theorem 1, this is more interesting to contemplate. We
discuss this issue and provide examples that constitute our present
understanding of it.\\

\textbf{Acknowledgements. }In the last section, we make use of two
standard propositions, Prop. 7 and Prop. 8, not new to this paper.
We provide their proofs in the appendix. We thank Asaf Katz and Yuval
Peres for sketching to the author the proof of Prop. 7, and Hillel
Furstenberg for showing him the proof of Prop. 8 as well as providing
other valuable remarks. In addition, we want to thank Mike Hochman
for generously sharing with the author his example appearing in the
second subsection of the appendix. Section 5 relies
on it.

\section{The Proofs of Theorem 1 and Theorem 2}

Theorem 1 follows right away from Propositions 3 and 4. Proposition
3 is (ii)$\rightarrow$(iii) and Proposition 4 is (i)$\rightarrow$(ii).
(iii)$\rightarrow$(i) is trivial as (i) is a special case of (iii)
(by taking the stationary process with values which are the complex
conjugate of the values of the stationary process defined by $\mu$).\\
\\

\textbf{Proposition 3:} \textit{Let $x=\left(x_{n}\right)_{n\in\mathbb{N}}$
be a complex sequence satisfying \\\\ $\sup_{N\in \mathbb{N}} \frac{\left|x_{1}\right|^2+\dots+\left|x_{N}\right|^2}{N}<\infty$.
Assume that for any $\varepsilon>0$ there exists $N'>0$ such that
for any $N''>N'$ there exists $T'>0$ for which}

\textit{$\frac{\left|\left\{ 0<\tau\leq T\,:\,\exists N'\leq N\le N''\,\left|\frac{1}{N}\stackrel[n=1]{N}{\sum}x_{n+\tau}\overline{x_{n}}\right|\geq\varepsilon\right\} \right|}{T}<\varepsilon$
for every $T>T'$ (we denote this condition by $\left(*\right)$).
Then $x$ is a pointwise cancellation sequence.}\\

\textbf{Proof:} Assuming the theorem is false, there is a $\sigma$-invariant
Borel probability measure on $\mathbb{C}^{\mathbb{N}}$ for which $\limsup_{N\rightarrow\infty}\left|\frac{1}{N}\stackrel[n=1]{N}{\sum}x_{n}y_{n}\right|>0$
on a set of positive measure. By considering its ergodic decomposition, there exists an
ergodic $\sigma$-invariant Borel probability measure $\nu$ on $\mathbb{C}^{\mathbb{N}}$,
a positive $\nu$-measure set $B'\subseteq\mathbb{C}^{\mathbb{N}}$
and some $a>0$ such that \\ $\limsup_{N\rightarrow\infty}\left|\frac{1}{N}\stackrel[n=1]{N}{\sum}x_{n}y_{n}\right|>a$
for all $y\in B'$. We may also assume that $\left\Vert \pi_{1}\right\Vert _{2}=1$
in $L_{\nu}^{2}\left(\mathbb{C}^{\mathbb{N}}\right)$ (by performing
dilation).\\

There is a positive $\nu$-measure set $B\subseteq B'$ and a sequence
of intervals $\left(L_{j},M_{j}\right)$ satisfying
\[
0<L_{1}<M_{1}<L_{2}<M_{2}<\dots
\]
such that for every $y\in B$ and interval $j\in\mathbb{N}$ there
exists some \\\\ $n_{j}\left(y\right)\in\left(L_{j},M_{j}\right)$ for
which $\left|\frac{1}{n_{j}\left(y\right)}\stackrel[n=1]{n_{j}\left(y\right)}{\sum}x_{n}y_{n}\right|>a$.
\\

We now fix $\delta=\frac{1}{3}$. There exists, by ergodicity, some
$K>0$ such that 
\[
\nu\left(\stackrel[k=0]{K-1}{\cup}\sigma^{-k}\left(B\right)\right)>1-\frac{\delta}{2}.
\]
We denote by $G$ the set of all $y\in\mathbb{C}^{\mathbb{N}}$ satisfying
\[
\frac{1}{N}\stackrel[n=0]{N-1}{\sum}1_{\stackrel[k=0]{K-1}{\cup}\sigma^{-k}\left(B\right)}\left(\sigma^{n}\left(y\right)\right)\xrightarrow[N\rightarrow\infty]{}\nu\left(\stackrel[k=0]{K-1}{\cup}\sigma^{-k}\left(B\right)\right)
\]
 and $\frac{1}{N}\stackrel[n=1]{N}{\sum}\left|y_{n}\right|^{2}\xrightarrow[N\rightarrow\infty]{}\left\Vert \pi_{1}\right\Vert _{2}^{2}=1$.
$\nu\left(G\right)=1$ so in particular $G$ is non-empty and this
is all we shall need. \\
\\

\uline{Our Goal:} Given any $y\in G$ and $J\in\mathbb{N}$, for
a large enough $D$ (depending on $y$ and $J$) we aim to construct
a sequence \[c^{\left(1\right)}\left(y\right),c^{\left(2\right)}\left(y\right),\dots,c^{\left(J\right)}\left(y\right)\in\mathbb{C}^{D} \]
satisfying the following (when considering these vectors with the
inner product of the normalized uniform measure on $\left\{ 1,\dots,D\right\} $):
\begin{enumerate}
\item $\left\Vert c^{\left(j\right)}\left(y\right)\right\Vert _{2}^{2}\leq\sup _{N\in\mathbb{N}}\frac{\left|x_{1}\right|^{2}+\dots+\left|x_{N}\right|^{2}}{N}$
for every $j$.
\item $\left\langle c^{\left(j_{1}\right)}\left(y\right),c^{\left(j_{2}\right)}\left(y\right)\right\rangle \leq\frac{\delta}{J}$
for every $j_{1}\neq j_{2}$.
\item $\left\langle c^{\left(j\right)}\left(y\right),\overline{\left(y_{n}\right)_{n=1}^{D}}\right\rangle >\left(1-\delta\right)a$
for every $j$.\\
\end{enumerate}
This implies that $c\left(y\right)=c^{\left(1\right)}\left(y\right)+\dots+c^{\left(J\right)}\left(y\right)$
satisfies on the one hand
\[
\left\Vert c\left(y\right)\right\Vert _{2}^{2}=\left\Vert c^{\left(1\right)}\left(y\right)+\dots+c^{\left(J\right)}\left(y\right)\right\Vert _{2}^{2}\leq\sup_{N\in\mathbb{N}}\frac{\left|x_{1}\right|^{2}+\dots+\left|x_{N}\right|^{2}}{N}J+\delta\left(J-1\right),
\]

And on the other hand, for $D$ large enough
\[
\left(1-\delta\right)aJ<\left\langle c^{\left(1\right)}\left(y\right)+\dots+c^{\left(J\right)}\left(y\right),\overline{\left(y_{n}\right)_{n=1}^{D}}\right\rangle \leq\left\Vert c\left(y\right)\right\Vert _{2}\sqrt{\frac{1}{D}\stackrel[n=1]{D}{\sum}\left|y_{n}\right|^{2}}
\]
\[
<\left(1+\delta\right)\left\Vert \pi_{1}\right\Vert _{2}\left\Vert c\left(y\right)\right\Vert _{2}=\left(1+\delta\right)\left\Vert c\left(y\right)\right\Vert _{2}.
\]

Hence we obtain that $\left(1-\delta\right)^{2}a^{2}J^{2}<\left(1+\delta\right)^{2}\left(\sup_{N\in\mathbb{N}}\frac{\left|x_{1}\right|^{2}+\dots+\left|x_{N}\right|^{2}}{N} J+\delta\left(J-1\right)\right)$
for all $J\in\mathbb{N}$, and this leads to the contradiction $a=0$.\\
\\

\uline{Constructing the sequence \mbox{$c^{\left(1\right)}\left(y\right),c^{\left(2\right)}\left(y\right),\dots,c^{\left(J\right)}\left(y\right)\in\mathbb{C}^{D}$}}:
Fixing some $y\in G$ and $J\in\mathbb{N}$, we construct these vectors
for a large enough $D$. By deleting some of the intervals $\left(L_{j},M_{j}\right)$
we may assume without loss of generality that $\frac{K}{L_{1}}<\frac{\delta}{6J},\,\frac{K+M_{j}}{L_{j+1}}<\frac{\delta}{6J}$,
and also, placing $\varepsilon=\frac{\delta}{6J^{2}K}$ in $\left(*\right)$,
the following requirements that we shall denote by $\left(\vartriangle\right)$:

$L_{1}\geq N'$, for every $1\leq j<J$, we choose $N''$ to be $M_{j}$
and require that $L_{j+1}\geq M_{j}+T'$ ($N',N''$ and $T'$ are
as in $\left(*\right)$ - notice that $N''$ and $T'$ depend on $j$
) .\\

We define $c^{\left(j\right)}\left(y\right)$ recursively on $j$
from $J$ downwards to $1$. As the base case, $c^{\left(J\right)}\left(y\right)$
is defined as follows: we denote by $l_{1}^{\left(J\right)}\left(y\right)$
the minimal \\\\ $0\leq l\leq D-M_{J}$ for which $\sigma^{l}\left(y\right)\in B$
and define 
\[
c_{l_{1}^{\left(J\right)}\left(y\right)}^{\left(J\right)}\left(y\right),\dots,c_{l_{1}^{\left(J\right)}\left(y\right)+n_{J}\left(\sigma^{l_{1}^{\left(J\right)}\left(y\right)}\left(y\right)\right)-1}^{\left(J\right)}\left(y\right)
\]
 to be equal to 
\[
\left(\omega_{1}^{\left(J\right)}x_{1}\right),\dots,\left(\omega_{1}^{\left(J\right)}x_{n_{J}\left(\sigma^{l_{1}^{\left(J\right)}\left(y\right)}\left(y\right)\right)}\right)
\]
 for the $\left|\omega_{1}^{\left(J\right)}\right|=1$ that satisfies
that 
\[
\frac{\omega_{1}^{\left(J\right)}}{n_{J}\left(\sigma^{l_{1}^{\left(J\right)}\left(y\right)}\left(y\right)\right)}\stackrel[n=1]{n_{J}\left(\sigma^{l_{1}^{\left(J\right)}\left(y\right)}\left(y\right)\right)}{\sum}x_{n}\cdot\left(\sigma^{l_{1}^{\left(J\right)}\left(y\right)}\left(y\right)\right)_{n}
\]
 is a positive number (greater than $a$). Then denote by $l_{2}^{\left(J\right)}\left(y\right)$
the minimal $l_{1}^{\left(J\right)}\left(y\right)+n_{J}\left(\sigma^{l_{1}^{\left(J\right)}\left(y\right)}\left(y\right)\right)\leq l\leq D-M_{J}$
for which $\sigma^{l}\left(y\right)\in B$ and define 
\[
c_{l_{2}^{\left(J\right)}\left(y\right)}^{\left(J\right)}\left(y\right),\dots,c_{l_{2}^{\left(J\right)}\left(y\right)+n_{J}\left(\sigma^{l_{2}^{\left(J\right)}\left(y\right)}\left(y\right)\right)-1}^{\left(J\right)}\left(y\right)
\]
 to be equal to 
\[
\left(\omega_{2}^{\left(J\right)}x_{1}\right),\dots,\left(\omega_{2}^{\left(J\right)}x_{n_{J}\left(\sigma^{l_{2}^{\left(J\right)}\left(y\right)}\left(y\right)\right)}\right)
\]
 for the $\left|\omega_{2}^{\left(J\right)}\right|=1$ satisfying
that 
\[
\frac{\omega_{2}^{\left(J\right)}}{n_{J}\left(\sigma^{l_{2}^{\left(J\right)}\left(y\right)}\left(y\right)\right)}\stackrel[n=1]{n_{J}\left(\sigma^{l_{2}^{\left(J\right)}\left(y\right)}\left(y\right)\right)}{\sum}x_{n}\cdot\left(\sigma^{l_{2}^{\left(J\right)}\left(y\right)}\left(y\right)\right)_{n}
\]
 is a positive number (greater than $a$). Thus continuing until reaching
step $R_{J}$ for which after $l_{R_{J}}^{\left(J\right)}\left(y\right)+n_{J}\left(\sigma^{l_{R_{J}}^{\left(J\right)}\left(y\right)}\left(y\right)\right)-1$
there is no such $l$. In all the rest of the indices of $c^{\left(J\right)}\left(y\right)$
we set the value to be zero.\\

For indices $1\leq n\leq D$ that do not belong to the union \\\\ $\stackrel[i=1]{R_{J}}{\cup}\left\{ l_{i}^{\left(J\right)}\left(y\right)-K+1,\dots,l_{i}^{\left(J\right)}\left(y\right),\dots,l_{i}^{\left(J\right)}\left(y\right)+n_{J}\left(\sigma^{l_{i}^{\left(J\right)}\left(y\right)}\left(y\right)\right)-1\right\} $,
necessarily $\sigma^{n}\left(y\right)\notin\stackrel[k=0]{K-1}{\cup}\sigma^{-k}\left(B\right)$
or $n\in\left\{ D-M_{J}-K+1,\dots,D\right\} $. Thus, for a large
$D$, there are more than $\left(1-\frac{\delta}{2}\right)D$ indices
$1\leq n\leq D$ that belong to that union, and thus there are more
than 
\[
\left(1-\frac{K}{L_{1}}\right)\left(1-\frac{\delta}{2}\right)D>\left(1-\frac{K}{L_{1}}-\frac{\delta}{2}\right)D>\left(1-\frac{\delta}{6J}-\frac{\delta}{2}\right)D>\left(1-\frac{\delta}{2J}-\frac{\delta}{2}\right)D
\]
 indices $1\leq n\leq D$ that belong to the union $\stackrel[i=1]{R_{J}}{\cup}\left\{ l_{i}^{\left(J\right)}\left(y\right),\dots,l_{i}^{\left(J\right)}\left(y\right)+n_{J}\left(\sigma^{l_{i}^{\left(J\right)}\left(y\right)}\left(y\right)\right)-1\right\} $
and also satisfy $\sigma^{n}\left(y\right)\in\stackrel[k=0]{K-1}{\cup}\sigma^{-k}\left(B\right)$.
Likewise, in the step corresponding to each $1\leq j_{0}\leq J$ in
the recursion we will prove inductively that there are more than $\left(1-\left(J-j_{0}+1\right)\frac{\delta}{2J}-\frac{\delta}{2}\right)D$
indices $1\leq n\leq D$ that belong to the union $\stackrel[i=1]{R_{j_{0}}}{\cup}\left\{ l_{i}^{\left(j_{0}\right)}\left(y\right),\dots,l_{i}^{\left(j_{0}\right)}\left(y\right)+n_{j_{0}}\left(\sigma^{l_{i}^{\left(j_{0}\right)}\left(y\right)}\left(y\right)\right)-1\right\} $
and also satisfy $\sigma^{n}\left(y\right)\in\stackrel[k=0]{K-1}{\cup}\sigma^{-k}\left(B\right)$.\\
\\
We turn to the step of recursion for $1\leq j_{0}<J$. $c^{\left(j_{0}\right)}\left(y\right)$
is defined with $n_{j_{0}}$ as we did in the base case with $c^{\left(J\right)}\left(y\right)$
and $n_{J}$, with an additional requirement when determining $l_{i}^{\left(j_{0}\right)}\left(y\right)$
that we now describe. We determine $l_{i}^{\left(j_{0}\right)}\left(y\right)$
similarly to how we determined $l_{i}^{\left(J\right)}\left(y\right)$,
only that we want \\ $1\leq l\leq D-M_{j_{0}}$ to not only satisfy $\sigma^{l}\left(y\right)\in B$,
but to also satisfy that $\left\{ l,\dots,l+n_{j_{0}}\left(\sigma^{l}\left(y\right)\right)-1\right\} \subseteq\left\{ l_{i}^{\left(j_{0}+1\right)}\left(y\right),\dots,l_{i}^{\left(j_{0}+1\right)}\left(y\right)+n_{j_{0}+1}\left(\sigma^{l_{i}^{\left(j_{0}+1\right)}\left(y\right)}\left(y\right)\right)-1\right\} $
for some $i$, and that for all $j_{0}<j\leq J$
\[
\frac{1}{n_{j_{0}}\left(\sigma^{l}\left(y\right)\right)}\left|x_{1}c_{l}^{\left(j\right)}\left(y\right)+\dots+x_{n_{j_{0}}\left(\sigma^{l}\left(y\right)\right)}c_{l+n_{j_{0}}\left(\sigma^{l}\left(y\right)\right)-1}^{\left(j\right)}\left(y\right)\right|<\frac{\delta}{2J}.
\]
\\

As before, the rest of the indices are set to zero.\\

If an index $1\leq n\leq D$ does not belong to either \\\\ $\stackrel[i=1]{R_{j_{0}}}{\cup}\left\{ l_{i}^{\left(j_{0}\right)}\left(y\right)-K+1,\dots,l_{i}^{\left(j_{0}\right)}\left(y\right),\dots,l_{i}^{\left(j_{0}\right)}\left(y\right)+n_{j_{0}}\left(\sigma^{l_{i}^{\left(j_{0}\right)}\left(y\right)}\left(y\right)\right)-1\right\} $
or \\\\ $\stackrel[k=0]{K-1}{\cup}\sigma^{-k}\left(B\right)$ then it must
belong to at least one of the following sets:\\

\[
A_{1}=\left\{ 1\leq n\leq D\,:\,\sigma^{n}\left(y\right)\notin\stackrel[k=0]{K-1}{\cup}\sigma^{-k}\left(B\right)\right.
\]
\[
\left.or\,\,n\notin\stackrel[i=1]{R_{j_{0}+1}}{\cup}\left\{ l_{i}^{\left(j_{0}+1\right)}\left(y\right),\dots,l_{i}^{\left(j_{0}+1\right)}\left(y\right)+n_{j_{0}}\left(\sigma^{l_{i}^{\left(j_{0}+1\right)}\left(y\right)}\left(y\right)\right)-1\right\} \right\} ,
\]

\[
A_{2}=\left\{ 1\leq n\leq D\,:\,n\notin A_{1},\,\forall i\,\,\left\{ n,\dots,n+K+M_{j_{0}}-1\right\} \right.
\]
\[
\left.\nsubseteq\left\{ l_{i}^{\left(j_{0}+1\right)}\left(y\right),\dots,l_{i}^{\left(j_{0}+1\right)}\left(y\right)+n_{j_{0}+1}\left(\sigma^{l_{i}^{\left(j_{0}+1\right)}\left(y\right)}\left(y\right)\right)-1\right\} \right\} ,
\]

\[
A_{3}=\stackrel[k=0]{K-1}{\cup}B_{k}\,\,\,\,where\,\,\,\,B_{k}=\left\{ 1\leq n\leq D\,:\,\text{\ensuremath{\,\sigma^{n+k}\left(y\right)\in B,}}\right.
\]

\[
\,\exists i\,\,\left\{ n+k,\dots,n+k+n_{j_{0}}\left(\sigma^{n+k}\left(y\right)\right)-1\right\} \subseteq\left\{ l_{i}^{\left(j_{0}+1\right)}\left(y\right),\dots,l_{i}^{\left(j_{0}+1\right)}\left(y\right)+n_{j_{0}+1}\left(\sigma^{l_{i}^{\left(j_{0}+1\right)}\left(y\right)}\left(y\right)\right)-1\right\} ,
\]

\[
\left.\,\exists j_{_{0}}<j\leq J\,\,\frac{1}{n_{j_{0}}\left(\sigma^{l}\left(y\right)\right)}\left|x_{1}c_{n+k}^{\left(j\right)}\left(y\right)+\dots+x_{n_{j_{0}}\left(\sigma^{l}\left(y\right)\right)}c_{n+k+n_{j_{0}}\left(\sigma^{n+k}\left(y\right)\right)-1}^{\left(j\right)}\left(y\right)\right|\geq\frac{\delta}{2J}\right\} .
\]
\\

For a large $D$, by the induction hypothesis $\left|A_{1}\right|<\left(\left(J-j_{0}\right)\frac{\delta}{2J}+\frac{\delta}{2}\right)D$.\\

As for $\left|A_{2}\right|$, since $n\notin A_{1}$, there exists
some $i$ for which \\\\ $n\in\left\{ l_{i}^{\left(j_{0}+1\right)}\left(y\right),\dots,l_{i}^{\left(j_{0}+1\right)}\left(y\right)+n_{j_{0}+1}\left(\sigma^{l_{i}^{\left(j_{0}+1\right)}\left(y\right)}\left(y\right)\right)-1\right\} $,
however

\[
\left\{ n,\dots,n+K+M_{j_{0}}-1\right\} \nsubseteq\left\{ l_{i}^{\left(j_{0}+1\right)}\left(y\right),\dots,l_{i}^{\left(j_{0}+1\right)}\left(y\right)+n_{j_{0}+1}\left(\sigma^{l_{i}^{\left(j_{0}+1\right)}\left(y\right)}\left(y\right)\right)-1\right\} .
\]
 Thus $n$ is one of the largest $K+M_{J_{0}}-1$ numbers in 
\[
\left\{ l_{i}^{\left(j_{0}+1\right)}\left(y\right),\dots,l_{i}^{\left(j_{0}+1\right)}\left(y\right)+n_{j_{0}+1}\left(\sigma^{l_{i}^{\left(j_{0}+1\right)}\left(y\right)}\left(y\right)\right)-1\right\} ,
\]
 and $\left|A_{2}\right|$ is bounded from above by $\frac{K+M_{j_{0}}}{L_{j_{0}+1}}D$.\\

Regarding $\left|A_{3}\right|$, let us first define for every $j_{0}<j\leq J$
the set $B_{k,j}$ as we defined $B_{k}$ except without the expression ``$\exists j_{_{0}}<j\leq J$''. So $B_{k}=\stackrel[j_{0}+1]{J}{\cup}B_{k,j}$.
By $\left(\triangle\right)$, from those $n$ that satisfy the first
two requirements in the definition of $B_{k,j}$ only a fraction of
less than $\varepsilon$ satisfies the third, thus certainly \\$\left|B_{k,j}\right|<\varepsilon D$.
Hence $\left|B_{k}\right|<J\varepsilon D=J\cdot\frac{\delta}{6J^{2}K}\cdot D=\frac{\delta}{6JK}D$
and thus \\ $\left|A_{3}\right|<\frac{\delta}{6J}D<\frac{\delta}{2J}D$.\\

Overall, we deduce that the number of indices $1\leq n\leq D$ that
do not belong to either

$\stackrel[i=1]{R_{j_{0}}}{\cup}\left\{ l_{i}^{\left(j_{0}\right)}\left(y\right)-K+1,\dots,l_{i}^{\left(j_{0}\right)}\left(y\right),\dots,l_{i}^{\left(j_{0}\right)}\left(y\right)+n_{j_{0}}\left(\sigma^{l_{i}^{\left(j_{0}\right)}\left(y\right)}\left(y\right)\right)-1\right\} $
or $\stackrel[k=0]{K-1}{\cup}\sigma^{-k}\left(B\right)$ is not greater
than
\[
\left|A_{1}\right|+\left|A_{2}\right|+\left|A_{3}\right|<\left(\left(J-j_{0}\right)\frac{\delta}{2J}+\frac{\delta}{2}\right)D+\frac{K+M_{j_{0}}}{L_{j_{0}+1}}D+\frac{\delta}{6J}D
\]
\[
<\left(\left(J-j_{0}\right)\frac{\delta}{2J}+\frac{\delta}{2}+\frac{\delta}{6J}+\frac{\delta}{6J}\right)D.
\]
And since $\frac{K}{L_{j_{0}}}<\frac{K}{L_{1}}<\frac{\delta}{6J}$,
the number of indices $1\leq n\leq D$ that do not belong to either

$\stackrel[i=1]{R_{j_{0}}}{\cup}\left\{ l_{i}^{\left(j_{0}\right)}\left(y\right),\dots,l_{i}^{\left(j_{0}\right)}\left(y\right)+n_{j_{0}}\left(\sigma^{l_{i}^{\left(j_{0}\right)}\left(y\right)}\left(y\right)\right)-1\right\} $
or $\stackrel[k=0]{K-1}{\cup}\sigma^{-k}\left(B\right)$ is smaller
than
\[
\frac{\delta}{6J}D+\left(\left(J-j_{0}\right)\frac{\delta}{2J}+\frac{\delta}{2}+\frac{\delta}{6J}+\frac{\delta}{6J}\right)D
\]
\[
=\left(\left(J-j_{0}+1\right)\frac{\delta}{2J}+\frac{\delta}{2}\right)D.
\]
\\

Up to now, we have constructed $c^{\left(1\right)}\left(y\right),\dots,c^{\left(J\right)}\left(y\right)$.
By the construction it is clear that $\left\Vert c^{\left(j\right)}\left(y\right)\right\Vert _{2}^{2}\leq\sup _{N\in\mathbb{N}}\frac{\left|x_{1}\right|^{2}+\dots+\left|x_{N}\right|^{2}}{N}$
for every $j$ and that $\left\langle c^{\left(j_{1}\right)}\left(y\right),c^{\left(j_{2}\right)}\left(y\right)\right\rangle \leq\frac{\delta}{J}$
for every $j_{1}\neq j_{2}$. We have also proved that for every $j$
there are more than $\left(1-\left(J-j_{0}+1\right)\frac{\delta}{2J}-\frac{\delta}{2}\right)D$
indices $1\leq n\leq D$ that belong to $\stackrel[i=1]{R_{j}}{\cup}\left\{ l_{i}^{\left(j\right)}\left(y\right),\dots,l_{i}^{\left(j\right)}\left(y\right)+n_{j}\left(\sigma^{l_{i}^{\left(j\right)}\left(y\right)}\left(y\right)\right)-1\right\} $.
With the aid of this fact we can infer that $\left\langle c^{\left(j\right)}\left(y\right),\overline{\left(y_{n}\right)_{n=1}^{D}}\right\rangle >\left(1-\delta\right)a$:\\

\[
\left|\frac{1}{D}\stackrel[n=1]{D}{\sum}c_{n}^{\left(j\right)}\left(y\right)y_{n}\right|=\left|\frac{n_{j}\left(\sigma^{l_{1}^{\left(j\right)}\left(y\right)}\left(y\right)\right)}{D}\cdot\frac{1}{n_{j}\left(\sigma^{l_{1}\left(y\right)}\left(y\right)\right)}\stackrel[n=1]{n_{j}\left(\sigma^{l_{1}^{\left(j\right)}\left(y\right)}\left(y\right)\right)}{\sum}x_{n}\cdot\left(\sigma^{l_{1}^{\left(j\right)}\left(y\right)}\left(y\right)\right)_{n}\right|+\dots
\]
\[
+\left|\frac{n_{j}\left(\sigma^{l_{R_{j}}^{\left(j\right)}}\left(y\right)\right)}{D}\cdot\frac{1}{n_{j}\left(\sigma^{l_{R_{j}}^{\left(j\right)}}\left(y\right)\right)}\stackrel[n=1]{n_{J}\left(\sigma^{l_{R_{j}}^{\left(j\right)}\left(y\right)}\left(y\right)\right)}{\sum}x_{n}\cdot\left(\sigma^{l_{R_{j}}^{\left(j\right)}\left(y\right)}\left(y\right)\right)_{n}\right|
\]

\[
>\left(\frac{n_{j}\left(\sigma^{l_{1}^{\left(j\right)}\left(y\right)}\left(y\right)\right)}{D}+\dots+\frac{n_{j}\left(\sigma^{l_{R_{j}}^{\left(j\right)}\left(y\right)}\left(y\right)\right)}{D}\right)a>\left(1-\left(J-j_{0}+1\right)\frac{\delta}{2J}-\frac{\delta}{2}\right)a\geq\left(1-\delta\right)a.\,\,\blacksquare
\]
\\
\\\\
\textbf{Proposition 4:} \textit{Let $x=\left(x_{n}\right)_{n\in\mathbb{N}}$
be a complex sequence generic relative to a $\sigma$-invariant Borel
probability measure $\mu$ on $\mathbb{C}^{\mathbb{N}}$ that satisfies \\\\
$\frac{1}{T}\stackrel[n=1]{T-1}{\sum}x_{n}\overline{\xi_{n}}\xrightarrow[T\rightarrow\infty]{}0$
for $\xi\in\mathbb{C}^{\mathbb{N}}$ $\mu$-almost-surely. Then for
any $\varepsilon>0$ there exists $N'>0$ such that for any $N''>N'$
there exists $T'>0$ for which $\frac{\left|\left\{ 0<\tau\leq T\,:\,\exists N'\leq N\le N''\,\left|\frac{1}{N}\stackrel[n=1]{N}{\sum}x_{n+\tau}\overline{x_{n}}\right|\geq\varepsilon\right\} \right|}{T}<\varepsilon$
for every $T>T'$.}\\

\textbf{Proof:} Given an $\varepsilon>0$ and $\tilde{N}\in\mathbb{N}$,
we denote by $A_{\varepsilon}^{\tilde{N}}\subseteq\mathbb{C^{\mathbb{N}}}$
the measurable set containing all points $\xi\in\mathbb{C}^{\mathbb{N}}$
for which $\left|\frac{1}{N}\stackrel[n=1]{N-1}{\sum}x_{n}\overline{\xi_{n}}\right|<\varepsilon$
for all $N\geq\tilde{N}.$ For a large enough $\tilde{N}$ we have
$\mu\left(A_{\varepsilon}^{\tilde{N}}\right)>1-\varepsilon$. We denote
the minimal such $\tilde{N}$ by $N'$. For every $N''\geq N'$, the
set \\\\ $B_{\varepsilon}^{N',N''}=\left\{ \xi\in\mathbb{C}^{\mathbb{N}}\,:\,\forall N'\leq N\le N''\,\left|\frac{1}{N}\stackrel[n=1]{N-1}{\sum}x_{n}\overline{\xi_{n}}\right|<\varepsilon\right\} $
is of $\mu$-measure greater than $1-\varepsilon$. Its complement
is a closed set, and hence, since $x$ is generic, there exists some
$T'\in\mathbb{N}$ such that $\frac{\left|\left\{ 0<\tau\leq T\,:\,\sigma^{\tau}\left(x\right)\notin B_{\varepsilon}^{N',N''}\right\} \right|}{T}<\varepsilon$
for every $T>T'$. But
\[
\left\{ 0<\tau\leq T\,:\,\sigma^{\tau}\left(x\right)\notin B_{\varepsilon}^{N',N''}\right\} =\left\{ 0<\tau\leq T\,:\,\exists N'\leq N\le N''\,\left|\frac{1}{N}\stackrel[n=1]{N}{\sum}x_{n+\tau}\overline{x_{n}}\right|\geq\varepsilon\right\} 
\]
for every $T$, and hence we are done. \textbf{$\blacksquare$}\\
\\

\section{Mean Cancellation Sequences}

By a similar proof to that of Theorem 1, the following proposition
can be obtained. \\

\textbf{Proposition 5:} \textit{Let $x=\left(x_{n}\right)_{n\in\mathbb{N}}$
be a complex sequence, and let $T_{k}$ be an increasing sequence
of integers for which $\sup_{k\in\mathbb{N}} \frac{\left|x_{1}\right|^{2}+\dots+\left|x_{T_{k}}\right|^{2}}{T_{k}} <\infty$
and $\lim_{k\rightarrow\infty}\frac{1}{T_{k}}\stackrel[n=1]{T_{k}}{\sum}\delta_{\sigma^{n}\left(x\right)}$
converges weak-{*} to some probability measure $\mu$. The following conditions
are equivalent:}

\textbf{\textit{(i) }}\textit{$\frac{1}{T_{k}}\stackrel[n=1]{T_{k}-1}{\sum}x_{n}\overline{\xi_{n}}\xrightarrow[k\rightarrow\infty]{}0$
for $\xi\in\mathbb{C}^{\mathbb{N}}$ $\mu$-almost-surely.}

\textbf{\textit{(ii)}}\textit{ For any $\varepsilon>0$ there exists
$s'>0$ such that for any $s''>s'$ there exists $k'>0$ for which }

\textit{$\frac{\left|\left\{ 0<\tau\leq T_{k}\,:\,\exists s'\leq s\le s''\,\left|\frac{1}{T_{s}}\stackrel[n=1]{T_{s}}{\sum}x_{n+\tau}\overline{x_{n}}\right|\geq\varepsilon\right\} \right|}{T_{k}}<\varepsilon$
for every $k>k'$.}

\textbf{\textit{(iii)}}\textit{ For any complex-valued stationary
process $\left(Y_{n}\left(\omega\right)\right)_{n\in\mathbb{Z}}$
of finite variance $\frac{1}{T_{k}}\stackrel[n=1]{T_{k}}{\sum}x_{n}Y_{n}\left(\omega\right)\xrightarrow[k\rightarrow\infty]{}0$} almost-surely.\\

\textbf{Proof of Proposition 5:} The proof is carried out similarly
to that of Theorem 1. In the proof of Prop. 3 the index $n_{j}\left(y\right)$
should be equal to a $T_{k}$. $\blacksquare$\\

The condition $\sup_{N\in \mathbb{N}} \frac{\left|x_{1}\right|^2+\dots+\left|x_{N}\right|^2}{N}<\infty$
implies \\ $\sup_{N\in \mathbb{N}} \frac{\left|x_{1}\right|+\dots+\left|x_{N}\right|}{N}<\infty$,
and the latter condition suffices to ensure that for every $\varepsilon>0$
there exists $M>0$ for which the measure according to $\frac{1}{T}\stackrel[n=1]{T}{\sum}\delta_{\sigma^{n}\left(x\right)}$
of the set $\left\{ z\in\mathbb{C}^{\mathbb{N}}\,:\,\left|z_{1}\right|\leq M\right\} $
is greater than $1-\varepsilon$ for all $T\in\mathbb{N}$. Thus,
in the weak-{*} topology, any partial limit of the sequence $\frac{1}{T}\stackrel[n=1]{T}{\sum}\delta_{\sigma^{n}\left(x\right)}$ is a probability measure (there is no escape of mass).\\

Hence Cor. 6 is an immediate conclusion from Prop. 5. It does not
presuppose any genericity of $x$.\\

\textbf{Corollary 6:} \textit{Let $x=\left(x_{n}\right)_{n\in\mathbb{N}}$
be a complex sequence satisfying \\ $\sup_{N\in \mathbb{N}} \frac{\left|x_{1}\right|^2+\dots+\left|x_{N}\right|^2}{N}<\infty$.
The following conditions are equivalent:}

\textbf{\textit{(i) }}\textit{For any increasing sequence of integers
$T_{k}$ there exists a sub-sequence $T_{k_{l}}$ that satisfies that
for any $\varepsilon>0$ there exists $r'>0$ such that for any $r''>r'$
there exists $l'>0$ for which} $\frac{\left|\left\{ 0<\tau\leq T_{k_{l}}\,:\,\exists r'\leq r\le r''\,\left|\frac{1}{T_{k_{r}}}\stackrel[n=1]{T_{k_{r}}}{\sum}x_{n+\tau}\overline{x_{n}}\right|\geq\varepsilon\right\} \right|}{T_{k_{l}}}<\varepsilon$
for every $l>l'$.

\textbf{\textit{(ii) }}\textit{For any increasing sequence of integers
$T_{k}$ there exists a sub-sequence $T_{k_{l}}$ that satisfies for
any complex-valued stationary process $\left(Y_{n}\left(\omega\right)\right)_{n\in\mathbb{Z}}$
of finite variance that $\frac{1}{T_{k_{l}}}\stackrel[n=1]{T_{k_{l}}}{\sum}x_{n}Y_{n}\left(\omega\right)\xrightarrow[l\rightarrow\infty]{}0$
almost-surely.}\\

We are ready to prove Theorem 2.\\

\textbf{Proof of Theorem 2: }Given a stationary process $\left(Y_{n}\left(\omega\right)\right)_{n\in\mathbb{Z}}$
of finite variance, applying the spectral theorem with the vector $Y_{1}$ and the unitary shift map, we obtain that $\left\Vert \frac{1}{T}\stackrel[n=1]{T}{\sum}x_{n}Y_{n}\left(\omega\right)\right\Vert _{2}\xrightarrow[T\rightarrow\infty]{}0$
if and only if \\ $\left\Vert \frac{1}{T}\stackrel[n=1]{T}{\sum}x_{n}z^{n}\right\Vert _{2}\xrightarrow[T\rightarrow\infty]{}0$
in $L^{2}$ on the unit circle equipped with the spectral measure.\\

\textbf{(i)$\rightarrow$(iv):} Assume $\frac{1}{T}\stackrel[n=1]{T}{\sum}x_{n}z^{n}\xrightarrow[T\rightarrow\infty]{}0$
for all $z\in\mathbb{C}$ on the unit circle. Because 
\[
\left|\frac{1}{T}\stackrel[n=1]{T}{\sum}x_{n}z^{n}\right|\leq\left|\frac{1}{T}\stackrel[n=1]{T}{\sum}\left|x_{n}\right|\right|,
\]

and $\sup_{T\in\mathbb{N}}\left|\frac{1}{T}\stackrel[n=1]{T}{\sum}\left|x_{n}\right|\right|<\infty$
, we can apply the dominated convergence theorem to deduce the desired
result.\\

\textbf{(iv)$\rightarrow$(iii):} A bounded random variable is of
finite variance.\\

\textbf{(iii)$\rightarrow$(i):} Let us assume the contrary, that
there exists some $z_{0}$ on the unit circle satisfying $\limsup_{T\rightarrow\infty}\left|\frac{1}{T}\stackrel[n=1]{T}{\sum}x_{n}z_{0}^{n}\right|>0$.
Considering $Y_{n}\left(\omega\right)$ to be $\omega z_{0}^{n}$
where the underlying probability space is the unit circle with the
uniform measure, the procedure above of applying the spectral theorem
yields a spectral measure possessing an atom in $z_{0}$, thus \\ $\left\Vert \frac{1}{T}\stackrel[n=1]{T}{\sum}x_{n}z^{n}\right\Vert _{2}\underset{T\rightarrow\infty}{\nrightarrow}0$
and thus  $\left\Vert \frac{1}{T}\stackrel[n=1]{T}{\sum}x_{n}Y_{n}\left(\omega\right)\right\Vert _{2}\underset{T\rightarrow\infty}{\nrightarrow}0$.\\

Up to now we proved (i)$\leftrightarrow$(iii)$\leftrightarrow$(iv),
and it is left to link these conditions to condition (ii) which is
condition (i) of Cor. 6. We will link them to condition (ii) of Cor.
6 (which is equivalent to the latter by that corollary).\\

\textbf{(condition (ii) of Cor. 6)$\rightarrow$ (iii):} Given a stationary
process $\left(Y_{n}\left(\omega\right)\right)_{n\in\mathbb{Z}}$
bounded by $M>0$,
\[
\left|\frac{1}{T}\stackrel[n=1]{T}{\sum}x_{n}Y_{n}\left(\omega\right)\right|\leq M\left|\frac{1}{T}\stackrel[n=1]{T}{\sum}\left|x_{n}\right|\right|,
\]
and thus, by the dominated convergence theorem, for any increasing
sequence of integers $T_{k}$ there exists a sub-sequence $T_{k_{l}}$ satisfying \\ $\left\Vert \frac{1}{T_{k_{l}}}\stackrel[n=1]{T_{k_{l}}}{\sum}x_{n}Y_{n}\left(\omega\right)\right\Vert _{2}\xrightarrow[l\rightarrow\infty]{}0$.
This implies that $\left\Vert \frac{1}{T}\stackrel[n=1]{T}{\sum}x_{n}Y_{n}\left(\omega\right)\right\Vert _{2}\xrightarrow[T\rightarrow\infty]{}0$.\\

\textbf{(iv)$\rightarrow$(condition (ii) of Cor. 6):} $\left\Vert \frac{1}{T}\stackrel[n=1]{T}{\sum}x_{n}Y_{n}\left(\omega\right)\right\Vert _{2}\xrightarrow[T\rightarrow\infty]{}0$
implies that $\left\Vert \frac{1}{T_{k}}\stackrel[n=1]{T_{k}}{\sum}x_{n}Y_{n}\left(\omega\right)\right\Vert _{2}\xrightarrow[k\rightarrow\infty]{}0$
for any increasing sequence of integers $T_{k}$. In particular, there
exists a sub-sequence $T_{k_{l}}$ for which the convergence is almost-surely.
$\blacksquare$\\

\section{Generic Cancellation Sequences and the Spectrum of $\pi_{1}$ Relative
to the Koopman Operator}

In \cite{key-3}, it is shown that for a $\sigma$-invariant Borel
probability measure $\mu$ on $\mathbb{C}^{\mathbb{N}}$ with $\pi_{1}$
square integrable, the spectrum of $\pi_{1}$ relative to the Koopman operator of $\sigma$
is continuous if and only if $\mu$-almost-every $x\in\mathbb{C}^{\mathbb{N}}$
satisfies $\frac{1}{T}\stackrel[n=0]{T-1}{\sum}x_{n}\overline{\xi_{n}}\xrightarrow[T\rightarrow\infty]{}0$
for $\xi\in\mathbb{C}^{\mathbb{N}}$   $\mu$-almost-surely. If $\mu$
is ergodic then $x\in\mathbb{C}^{\mathbb{N}}$ is $\mu$-almost-surely
generic relative to $\mu$, and hence when also these two conditions
hold then it is   $\mu$-almost-surely a pointwise cancellation sequence
(by Theorem 1'). However, we now give an example of a pointwise cancellation
sequence $x\in\mathbb{C}^{\mathbb{N}}$ which is generic relative
to a $\sigma$-invariant Borel probability measure $\mu$ but the
two equivalent conditions above do not hold for that $\mu$.\\

We make use of the following two standard propositions whose proofs
are found in the appendix.\\

\textbf{Proposition 7:} \textit{For every irrational $\beta$ and
real $\alpha\neq0$, the sequence $\left(e^{2\pi n\beta i},e^{2\pi\text{\ensuremath{\sqrt{n}}}\alpha i}\right)_{n=1}^{\infty}$
equidistributes on the 2-torus.}\\

\textbf{Proposition 8:} \textit{Let $\left(Y,\mathscr{B},T,\nu\right)$
be probability measure preserving system and let $\alpha$ be an irrational
number. If $h$ is an integrable function on the circle (relative
to the uniform measure) and $f\in L_{\nu}^{2}\left(Y\right)$ is orthogonal
in $L_{\nu}^{2}\left(Y\right)$ to all eigenfunctions (relative to
the Koopman operator) with eigenvalues of the form $e^{2\pi m\alpha i}$
for $0\neq m\in\mathbb{Z}$, then, under every joining of $\left(Y,\mathscr{B},T,\nu\right)$
and the system of the circle rotation by $e^{2\pi\alpha i}$, the
integral of $f\left(y\right)h\left(z\right)$ equals the product of
the integrals of each of these functions.}\\

\textbf{Example:} Let $\alpha$ be an irrational number, and we denote
$\omega=e^{2\pi\alpha i}$. The sequence $x=\text{\ensuremath{\left(e^{2\pi\text{\ensuremath{\left(n+\sqrt{n}\right)}}\alpha i}\right)}}_{n=1}^{\infty}$ equidistributes on the circle (by Prop. 7) and - because also  $\sqrt{n+1}-\sqrt{n}\xrightarrow[n\rightarrow\infty]{}0$
- considered as a point in $\mathbb{C}^{\mathbb{N}}$ it is generic
relative to the measure which is the push-forward of the uniform measure
of the circle through the map $z\mapsto\left(z\omega^{n}\right)_{n=1}^{\infty}$.
The system of this measure together with the shift is measure-theoretically
isomorphic to the system of the rotation of the circle by $\omega$,
and so it is not weak-mixing. We now prove that $x$ is a pointwise
cancellation sequence. Given an ergodic $\sigma$-invariant Borel
probability measure $\nu$ on $\mathbb{C}^{\mathbb{N}}$ for which
$\pi_{1}$ is square integrable, we can write in $L_{\nu}^{2}\left(\mathbb{C}^{\mathbb{N}}\right)$:
$\pi_{1}\left(y\right)=f\left(y\right)+g\left(y\right)$ where $g$
belongs to the closure of the span of the eigenfunctions of eigenvalues
$\omega^{m}$ for $0\neq m\in\mathbb{Z}$ and $f$ is orthogonal to
this span. It is straight-forward to prove using the Cauchy-Schwarz inequality that the functions  $h\in L_{\nu}^{2}\left(\mathbb{C}^{\mathbb{N}}\right)$ for which $\lim_{T\rightarrow\infty}\frac{1}{T}\stackrel[n=1]{T}{\sum}x_{n}h\left(\sigma^{n}\left(y\right)\right)=0$ $\nu$-almost-surely form a closed sub-space (this is done in \cite{key-2}),
hence to show that $\lim_{T\rightarrow\infty}\frac{1}{T}\stackrel[n=1]{T}{\sum}x_{n}g\left(\sigma^{n}\left(y\right)\right)=0$
we may consider without loss of generality that $g=g_{m}$ where $g_{m}$
is an eigenfunction of eigenvalue $\omega^{m}$.\\

The equality $g_{m}\left(\sigma^{n}\left(y\right)\right)=\omega^{nm}g_{m}\left(y\right)$
holds for $y\in\mathbb{C}^{\mathbb{N}}$ $\nu$-almost-surely. So
for $y\in\mathbb{C}^{\mathbb{N}}$ $\nu$-almost-surely
\[
\frac{1}{T}\stackrel[n=1]{T}{\sum}x_{n}y_{n}=\frac{1}{T}\stackrel[n=1]{T}{\sum}x_{n}f\left(\sigma^{n}\left(y\right)\right)+\frac{1}{T}\stackrel[n=1]{T}{\sum}x_{n}g_{m}\left(\sigma^{n}\left(y\right)\right)\xrightarrow[T\rightarrow\infty]{}0+0=0.
\]

The first limit vanishes for $y\in\mathbb{C}^{\mathbb{N}}$ $\nu$-almost-surely
by Prop. 8 and the fact that the expectation of the stationary process
defined by $x$ is $0$. Regarding the second limit, it is the limit
in $T$ of the expression 
\[
\frac{1}{T}\stackrel[n=1]{T}{\sum}x_{n}g_{m}\left(\sigma^{n}\left(y\right)\right)=\frac{1}{T}\stackrel[n=1]{T}{\sum}e^{2\pi\text{\ensuremath{\left(n+\sqrt{n}\right)}}\alpha i}\omega^{mn}g_{m}\left(y\right)=g_{m}\left(y\right)\frac{1}{T}\stackrel[n=1]{T}{\sum}e^{2\pi\text{\ensuremath{\left(n\left(m+1\right)+\sqrt{n}\right)}}\alpha i}.
\]

and $\frac{1}{T}\stackrel[n=1]{T}{\sum}e^{2\pi\text{\ensuremath{\left(m\left(n+1\right)+\sqrt{n}\right)}}\alpha i}\xrightarrow[N\rightarrow\infty]{}0$
by Prop. 7 (the case of $m=0$ is in fact the case $m=0$ in the proof
of Prop. 7). //\\
\\
For $x$ of the above example, it is interesting to look into the
behaviour of the expression 
\[
\frac{1}{N}\stackrel[n=1]{N}{\sum}x_{n+\tau}\overline{x_{n}}=\frac{1}{N}\stackrel[n=1]{N}{\sum}e^{2\pi\text{\ensuremath{\left(n+\tau+\sqrt{n+\tau}\right)}}\alpha i}e^{-2\pi\text{\ensuremath{\left(n+\sqrt{n}\right)}}\alpha i}=\frac{1}{T}\stackrel[n=1]{N}{\sum}e^{2\pi\text{\ensuremath{\left(\tau+\sqrt{n+\tau}-\sqrt{n}\right)}}\alpha i}.
\]

By Theorem 1, for any $\varepsilon>0$ there exists $N'>0$ such that
for any $N''>N'$ there exists $T'>0$ for which 
\[\frac{\left|\left\{ 0<\tau\leq T\,:\,\exists N'\leq N\le N''\,\left|\frac{1}{N}\stackrel[n=1]{N}{\sum}e^{2\pi\text{\ensuremath{\left(\tau+\sqrt{n+\tau}-\sqrt{n}\right)}}\alpha i}\right|\geq\varepsilon\right\} \right|}{T}<\varepsilon\]

for every $T>T'$. And yet,

$\frac{1}{N}\stackrel[n=1]{N}{\sum}e^{2\pi\text{\ensuremath{\left(\tau+\sqrt{n+\tau}-\sqrt{n}\right)}}\alpha i}\xrightarrow[N\rightarrow\infty]{}e^{2\pi\tau\alpha i}$
for all $\tau$.\\

An objection to this example can be made claiming that $x$ is not
in the support of the measure it is generic to and maybe if this additional
requirement would have been added then we could not find such an $x$.
We will show this is not the case with a new example that makes use
of the previous one.\\

\textbf{Example: }Consider any weak-mixing $\sigma$-invariant Borel
probability measure $\eta$ on $\mathbb{C}^{\mathbb{N}}$ of full
support with $\pi_{1}$ being square integrable (e.g. independent trials of any fully supported
distribution on $\mathbb{C}$ of finite variance). We take $u\in\mathbb{C}^{\mathbb{N}}$
to be any of its generic points that also satisfy $\frac{1}{T}\stackrel[n=0]{T-1}{\sum}u_{n}\xi_{n}\xrightarrow[T\rightarrow\infty]{}0$
for $\xi\in\mathbb{C}^{\mathbb{N}}$ $\eta$-almost-surely, and continue
with $x$ from the previous example. The stationary process defined
by $\eta$ is disjoint from the Kronecker stationary process that $x$ is generic to, and hence $\left(x,u\right)$
is generic relative to the product measure of the two processes. Thus
$x+u$ is generic relative to the ergodic process which is the sum
of the two processes taken independently, i.e. the convolution of the two measures (these are
measure on the additive group $\mathbb{C}^{\mathbb{N}}$). This new
measure is of full support and, according to it, $\pi_{1}$ is square
integrable and does not have a continuous spectrum relative to the
Koopman operator. We claim that for any complex-valued stationary
process $\left(Y_{n}\left(\omega\right)\right)_{n\in\mathbb{N}}$
of finite variance
\[\frac{1}{N}\stackrel[n=1]{N}{\sum}\left(x+u\right)_{n}Y_{n}\left(\omega\right)=\frac{1}{N}\stackrel[n=1]{N}{\sum}x_{n}Y_{n}\left(\omega\right)+\frac{1}{N}\stackrel[n=1]{N}{\sum}u_{n}Y_{n}\left(\omega\right)\xrightarrow[N\rightarrow\infty]{}0+0\]
almost-surely. The proof of the limit of the first summand was carried
out in the previous example, and the second one follows from Theorem
1'. //\\

This last example of $x+u$ fixed the problem of $x$ not belonging
to the support of the measure it is generic to, but, unlike $x$,
$x+u$ is not bounded. We do not have in our possession such an example
that is also bounded.\\

\section{A Sequence with Mean but not Pointwise Cancellation}

This section demonstrates the existence of weak-mixing bounded stationary
processes of zero expectation that admits a generic point which is
not a pointwise cancellation sequence. However, any generic point
of a bounded weak-mixing stationary process of zero expectation satisfies
condition (i) in Theorem 0 (by the disjointness of weak-mixing and Kronecker processes) and as such forms a mean cancellation
sequence.\\

The proof relies on a yet unpublished example of M. Hochman appearing
in the appendix, and this section assumes the reader is now familiar
with it. Hochman shows that there exists a weak-mixing system of invariant
measure $\mu$ which has a generic point that satisfies the following:
the points in the fiber belonging to it in the product of the system
with itself are $\mu$-almost-surely not generic relative to any $\sigma\times\sigma$-invariant
measure. In fact, he proves the existence of a family of such measures
which he refers to as weak-mixing measures\textit{ admitting tall
covers}.\\

We hereby just indicate how to modfiy what is done in Theorem 10 for
our current need.\\

We choose $\Lambda=\left\{ a,b\right\} $ ($a$ and $b$ being some
symbols) and take any weak-mixing $\sigma$-invariant Borel probability
measure $\mu$ admitting tall covers. Insetad of defining $x$ as
in Theorem 10, we define it in a simpler manner: $x_{i}=\left(a'_{n}\right)_{i}$
for $\left|a'_{n-1}\right|<i\leq\left|a'_{n}\right|$. Substitute
$a$ and $b$ with any two real numbers $a'<0<b'$ for which the corresponding
process $\mu'$ on $\left\{ a',b'\right\} ^{\mathbb{N}}$ is of zero
expectation and denote by $x'$ the point in $\left\{ a',b'\right\} ^{\mathbb{Z}}$
corresponding to $x$. Then in $\left(\mathbb{C}^{\mathbb{N}}\times\mathbb{C}^{\mathbb{N}},\sigma\times\sigma\right)$
the point $\left(x',\xi\right)$ has a sequence of orbital measure
that converges to the measure $\Delta_{\mu'}$ (as in the notation
of the proof of Theorem 10) for $\xi\in\mathbb{C}^{\mathbb{N}}$ $\mu'$-almost-surely
. Thus $\mu'$-almost-surely for $\xi\in\mathbb{C}^{\mathbb{N}}$
the sequence $\frac{1}{N}\stackrel[i=1]{N}{\sum}x'_{i}\overline{\xi}$
has a sub-sequence that converges to $\mathbb{E}_{\mu}\left(\pi_{1}^{2}\right)>0$.

\section{Appendix}

\subsection{Proofs of Prop. 7 and Prop. 8}

We present here the proofs of Prop. 7 and Prop. 8 from the last section.
These claims are standard, but appear here in order that this paper
be self-contained.\\

\textbf{Proof of Prop. 7:} Let $\beta$ be an irrational number and
$\alpha\neq0$ be a real number. It suffices to prove that for every
function $I:S^{1}\rightarrow\left\{ 0,1\right\} $ which is an indicator
of an arc of the circle between angles $2\pi a$ and $2\pi b$ for
$0\leq a<b\leq1$ and $m\in\mathbb{Z}$, the sequence $\left(e^{2\pi n\beta i},e^{2\pi\text{\ensuremath{\sqrt{n}}}\alpha i}\right)_{n=1}^{\infty}$
satisfies $\frac{1}{N}\stackrel[n=0]{N}{\sum}I\left(e^{2\pi\text{\ensuremath{\sqrt{n}}}\alpha i}\right)e^{2\pi mn\beta i}\xrightarrow[m\rightarrow\infty]{}\left\{ \begin{array}{ccc}
0 &  & m\neq0\\
b-a &  & m=0
\end{array}\right.$ (since step functions uniformly approximate continuous functions).\\

We first deal with the case $m=0$. Given an $\varepsilon>0$, we
fix a positive integer $d$ for which the distance between $d\alpha$
and the nearest integer is smaller than $\frac{\varepsilon}{2}$.
It suffices to prove that 

\[\left|\frac{1}{2Nd+d^{2}}\stackrel[n=N^{2}]{\left(N+d\right)^{2}-1}{\sum}I\left(e^{2\pi\text{\ensuremath{\sqrt{n}}}\alpha i}\right)-\left(b-a\right)\right|<\varepsilon\]

for every $N$ large enough. And this follows easily from the following
consideration: $\sqrt{N^{2}+k}=N+\frac{1}{2N}k+R\left(N,k\right)$,
$\left|R\left(N,k\right)\right|<\frac{1}{4N^{3}}$ for $k\geq0$ (by
Taylor's remainder theorem for the real variable $k$). So when $k$
ranges on the integers between $0$ and $2Nd+d^{2}-1$, the sequence
$\sqrt{N^{2}+k}\cdot\alpha$ traverses from $\sqrt{N^{2}}\cdot\alpha=N\alpha$
to $\sqrt{\left(N+d\right)^{2}}\cdot\alpha=\left(N+d\right)\alpha$
in nearly equal steps of $\frac{1}{2N}\alpha$ (for large $N$), and
$\left(N+d\right)\alpha-N\alpha=d\alpha$ is an integer up to a difference
of $\frac{\varepsilon}{2}$.\\

As for the case $m\neq0$, suppose $\varepsilon>0$ be given. It suffices
to prove that for a large enough $N$, if $e^{2\pi\text{\ensuremath{\sqrt{N}}}\alpha i}$
is inside the arc that $I$ indicates and $e^{2\pi\text{\ensuremath{\sqrt{N-1}}}\alpha i}$
is not, then $\left|\frac{1}{r\left(N\right)}\stackrel[n=N]{N+r\left(N\right)-1}{\sum}e^{2\pi mn\beta i}\right|<\varepsilon$
where $r\left(N\right)\in\mathbb{N}$ is the number for which $e^{2\pi\text{\ensuremath{\sqrt{N+k}}}\alpha i}$
stays inside the arc from $k=0$ precisely up to $k=r\left(N\right)-1$.
Evaluating a geometric series
\[
\left|\frac{1}{r\left(N\right)}\stackrel[n=N]{N+r-1}{\sum}e^{2\pi mn\beta i}\right|<\frac{1}{r\left(N\right)}\cdot\frac{1}{\left|1-e^{2\pi m\beta i}\right|}.
\]
 and observing that $r\left(N\right)\rightarrow \infty$ as $N\rightarrow\infty$
completes the proof. $\blacksquare$\\
\\

\textbf{Proof of Prop. 8:} Let $\lambda$ be a joining of the two
systems. We may assume without loss of generality that $h\left(z\right)=z^{m}$
for some $0\neq m\in\mathbb{Z}$. We need to prove that $\underset{Y\times S^{1}}{\int}f\left(y\right)z^{m}\,d\lambda\left(y,z\right)=0$.
By the invariance of $\lambda$, 
\[
\underset{Y\times S^{1}}{\int}f\left(y\right)z^{m}\,d\lambda\left(y,z\right)=\underset{Y\times S^{1}}{\int}\frac{1}{N}\stackrel[n=0]{N-1}{\sum}e^{2\pi mn\alpha i}z^{m}f\left(T^{n}\left(y\right)\right)\,d\lambda\left(y,z\right)
\]
\[
=\underset{Y\times S^{1}}{\int}z^{m}\frac{1}{N}\stackrel[n=0]{N-1}{\sum}e^{2\pi mn\alpha i}f\left(T^{n}\left(y\right)\right)\,d\lambda\left(y,z\right).
\]
On the other hand, $u\left(y\right)\mapsto e^{2\pi m\alpha i}f\left(T\left(y\right)\right)$
is an isometry from $L_{\nu}^{2}\left(Y\right)$ to itself, and thus
by the mean ergodic theorem $\lim_{N\rightarrow\infty}\frac{1}{N}\stackrel[n=0]{N-1}{\sum}e^{2\pi mn\alpha i}f\left(T^{n}\left(y\right)\right)$
exists in $L_{\nu}^{2}\left(Y\right)$ and equals the projection of
$f$ on the eigenspace - relative to the Koopman operator of $T$
- of the eigenvalue $e^{-2\pi m\alpha i}$. So by the assumption on
$f$ this yields that $\lim_{N\rightarrow\infty}\frac{1}{N}\stackrel[n=0]{N-1}{\sum}e^{2\pi mn\alpha i}f\left(T^{n}\left(y\right)\right)=0$
in $L_{\nu}^{2}\left(Y\right)$. Hence 
\[
\underset{Y\times S^{1}}{\int}f\left(y\right)z^{m}\,d\lambda\left(y,z\right)=\lim_{N\rightarrow\infty}\underset{Y\times S^{1}}{\int}z^{m}\frac{1}{N}\stackrel[n=0]{N-1}{\sum}e^{2\pi mn\alpha i}f\left(T^{n}\left(y\right)\right)\,d\lambda\left(y,z\right)=0.\,\,\blacksquare
\]
\\

\subsection{Hochman's Example}

Here we describe M. Hochman's example of a weak-mixing system of invariant
measure $\mu$ and its generic point whose fiber in the product of
the system with itself satisfies a peculiar property: the points on
it are $\mu$-almost-surely not generic relative to any $\sigma\times\sigma$-invariant
measure (the reader interested in this topic can also look into Section
3 in \cite{key-4}). It is actually not a construction of a specific
example, but rather proving the existence of a family of such invariant
measures.\\

In the following, $\Lambda$ is some finite alphabet. As usual, $\Lambda^{*}$
will denote the collection of all finite length words in $\Lambda$,
and if $\alpha\in\Lambda^{*}$ then $\left|\alpha\right|$ will denote
its length. Also, $\left[\alpha\right]$ is the set of points $y\in\Lambda^{\mathbb{N}}$
satisfying $y_{1}\dots y_{\left|\alpha\right|}=\alpha$.\\

\textbf{Definition:} Let $\mu$ be a $\sigma$-invariant Borel probability
measure on $\Lambda^{\mathbb{N}}$. For $\varepsilon>0,\,N\in\mathbb{N}$
we shall say that \textit{$\mu$ admits an $\left(\varepsilon,N\right)$-cover}
if there exists words $a_{1},\dots,a_{p}\in\Lambda^{*}$ satisfying
that $\left|a_{1}\right|>N$, $\left|a_{i+1}\right|>3\left|a_{i}\right|$
and $\mu\left(\stackrel[i=1]{p}{\cup}\left[a_{i}\right]\right)>1-\varepsilon$.
We shall say that \textit{$\mu$ admits tall covers} if it admits
an $\left(\varepsilon,N\right)$-cover for all $\left(\varepsilon,N\right)$.\\

Although measures admitting tall covers have zero entropy (this follows
easily from the definition of entropy on the obvious generating partition),
the family of those measures is a dense $G_{\delta}$-set relative
to the weak-{*} topology on the Borel probability measures on $\Lambda^{\mathbb{N}}$.
This is because, for every $\varepsilon$ and $N$, the set of measures
admitting $\left(\varepsilon,N\right)$-covers is open and contains
all periodic measures.\\

The family of Borel probability measures that both admit tall covers and are weak-mixing contains a dense $G_{\delta}$-set since also the set of all weak-mixing Borel probability measures on $\Lambda^{\mathbb{N}}$ contains a dense $G_{\delta}$-set (cf. the chapter ``Category'' in \cite{key-10}).\\

\textbf{Definition:} Let $\mu$ be a $\sigma$-invariant ergodic Borel
probability measure on $\Lambda^{\mathbb{N}}$, and let $\alpha\in\Lambda^{*}$.
For $\varepsilon>0$, we shall say that $\alpha$ is\textit{ $\varepsilon$-generic}
if each word $\beta$ of length not greater than $\frac{1}{\varepsilon}$
appears in it with a frequency closer than $\varepsilon$ to its $\mu$-probability,
i.e. 
\[
\left|\mu\left(\beta\right)-\frac{1}{\left|\alpha\right|}\left|\left\{ 1\leq i\leq\left|\alpha\right|-\left|\beta\right|+1\,:\,\alpha_{i}\dots\alpha_{i+\left|\beta\right|-1} = \beta \right\} \right|\right| < \varepsilon.
\]

For $M\in\mathbb{N}$, we shall say that $\alpha$ is \textit{strongly
$\left(\varepsilon,M\right)$-generic} if 
\[
\frac{1}{\left|\alpha\right|}\left|\left\{ 1\leq i<\left|\alpha\right|-M+1\,:\,\alpha_{i}\dots\alpha_{i+M-1}\,\,is\,\text{\ensuremath{\varepsilon-generic}}\right\} \right|>1-\varepsilon.
\]
\\
\textbf{Proposition 9:}\textit{ If $\mu$ is a $\sigma$-invariant
ergodic Borel probability measure on $\Lambda^{\mathbb{Z}}$ that
admits tall covers, then for every $\varepsilon>0$ and every large
enough $M$ there can be taken for any $N\in\mathbb{N}$ an $\left(\varepsilon,N\right)$-cover
comprised only of strongly $\left(\varepsilon,M\right)$-generic words.}\\

\textbf{Proof:} By the mean ergodic theorem, if $M$ is large enough,
there is a $\mu$-probability greater than $1-\varepsilon^{2}$ that
a point $y\in\Lambda^{\mathbb{N}}$ satisfies that $y_{1}\dots y_{M}$
is $\varepsilon$-generic. Given an $\left(\varepsilon,N\right)$-cover
$\tilde{E}$, we claim that if $N$ is large enough (which we may
suppose without loss of generality) then its elements that are not
strongly $\left(\varepsilon,M\right)$-generic occupy together a $\mu$-probability
of less than $3\varepsilon$. The mean ergodic theorem implies we
can find for every large enough $L$ some $y\in\Lambda^{\mathbb{N}}$
for which $y_{1}\dots y_{L}$ is both strongly $\left(\varepsilon,M\right)$-generic
and each of the elements of $\tilde{E}$ appears in it with a frequency
closer than $\frac{\varepsilon^{2}}{\left|\tilde{E}\right|}$ to its
$\mu$-probability.\\

Denote by $E$ the collection of elements of $\tilde{E}$ that are
not strongly $\left(\varepsilon,M\right)$-generic and assume without
loss of generality that none of which forms an initial segment of another.
Thus we can parse $y_{1}\dots y_{L}$ to read elemets of $E$ in the
following way: starting from $y_{1}$, we progress to the right until
we meet an element of $E$, then we keep progressing to the right
from the first index after that element ends until we meet another
element and so on. An index $1,\dots,L$ not covered by this parsing
cannot be a beginning of an element of $E$ in $y_{1}\dots y_{L}$ and
thus there are at most $\left(1-\mu\left(E\right)+\varepsilon^{2}\right)L$
such indices ($\mu\left(E\right)$ denotes, by abuse of notation,
the $\mu$-probability that the elements of $E$ occupy together).
Hence this parsing covers at least $\left(\mu\left(E\right)-\varepsilon^{2}\right)L$
of the indices $1,\dots,L$. For $N$ much larger than $M$ this means
that at least $\frac{1}{2}\left(\mu\left(E\right)-\varepsilon^{2}\right)L$
of the $M$-blocks in $y_{1}\dots y_{L}$ are not $\varepsilon$-generic.
So if $\mu\left(E\right)\geq3\varepsilon$ then $\frac{1}{2}\left(\mu\left(E\right)-\varepsilon^{2}\right)\geq\frac{3}{2}\varepsilon-\frac{1}{2}\varepsilon^{2}>\varepsilon$
which yields a contradiction to the fact that $y_{1}\dots y_{L}$
is strongly $\left(\varepsilon,M\right)$-generic.\\

Hence we can omitt from the cover the elements of $E$ and remain
with a $\left(4\varepsilon,N\right)$-cover comprised only of strongly
$\left(\varepsilon,M\right)$-generic words. $\blacksquare$\\

We shall also need an elementary lemma.\\

\textbf{Lemma 10:}\textit{ If $\alpha$ is strongly $\left(\varepsilon,M\right)$-generic
then it is $2\varepsilon$-generic.}\\

\textbf{Proof: }We parse $\alpha$ as follows. We find the minimal
$i$ for which $\alpha_{i}\dots\alpha_{i+M-1}$ is $\varepsilon$-generic,
include this word in the parsing, then progress to the least index greater than $i+M-1$ satisfying the same property, include that word and thus continue.
All words of length not greater than $\frac{1}{\varepsilon}$ appear
in frequency closer than $\varepsilon$ to their $\mu$-probability
among the indices that are included in the parsing as initial indices.
All the indices $i$ that are not included in this parsing are in
particular indices for which $\alpha_{i}\dots\alpha_{i+M-1}$ is not
$\varepsilon$-generic, hence there are less than $\varepsilon\left|\alpha\right|$
such indices. So they cannot distort the frequency of those words by
more than $\varepsilon$. $\blacksquare$\\

It remains to prove the following Theorem.\\

\textbf{Theorem 11:} \textit{Let $\mu$ be a $\sigma$-invariant ergodic
Borel probability measure on $\Lambda^{\mathbb{Z}}$ admitting tall
covers. Then there exists a generic point $x\in\Lambda^{\mathbb{N}}$
relative to $\mu$ such that the points composing $\left\{ x\right\} \times\Lambda^{\mathbb{N}}$
are $\mu$-almost-surely not generic relative to any $\sigma\times\sigma$-invariant
Borel probability measure.}\\

\textbf{Proof:} For $\varepsilon_{k}=\frac{1}{2^{k}}$ and a large
enough $M_{k}$ as in Prop. 9, we choose - inductively on $k\in\mathbb{N}$
- $\left(N_{k},\varepsilon_{k}\right)$-covers $a_{k,1},\dots,a_{k,p_{k}}$
with all elements strongly $\left(\varepsilon_{k},M_{k}\right)$-generic
and with $N_{k}$ greater than both $k\left|a_{k-1,p_{k-1}}\right|$
(for $k>1$) and $3\cdot2^{k}M_{k}$.\\

For notational convenience we denote the sequence made of concatenating
all the elements of these covers together in order $a_{1,1},\dots,a_{1,p_{1}},a_{2,1},\dots,a_{2,p_{2}},a_{3,1},\dots,a_{3,p_{3}},\dots$
by $\left\{ a'_{n}\right\} _{n=1}^{\infty}$. And by $k\left(n\right)$
we denote the $k$ corresponding to each $n$.\\

We are ready to construct $x$:
\[
x_{i}=\left\{ \begin{array}{ccc}
\left(a'_{n}\right)_{i} &  & \left|a'_{n-1}\right|<i\leq\left|a'_{n-1}\right|+\left\lceil \frac{1}{3}\left|a'_{n}\right|\right\rceil \\
\left(a'_{n}\right)_{i-1} &  & \left|a'_{n-1}\right|+\left\lceil \frac{1}{3}\left|a'_{n}\right|\right\rceil <i\leq\left|a'_{n}\right|
\end{array}\right..
\]
\\
Since $a'_{n}$ is strongly $\left(\frac{1}{2^{k\left(n\right)}},M_{k\left(n\right)}\right)$-generic,
every sub-block of length at least $\frac{1}{3}\left|a'_{n}\right|$
must at least be strongly $\left(\frac{4}{2^{k\left(n\right)}},M_{k\left(n\right)}\right)$-generic
(the percentage of ``good'' $M$-blocks in the sub-block need not
be less than $\frac{3}{2^{k\left(n\right)}}$ due to ``boundary effects'',
but the requirement $M_{k}<\frac{1}{2^{k}}\frac{N_{k}}{3}$ ensures
it is less than $\frac{4}{2^{k\left(n\right)}}$), so in particular
$\frac{8}{2^{k\left(n\right)}}$-generic by Lemma 10. This - together
with the fact that

$k\left(n\right)\xrightarrow[n\rightarrow\infty]{}\infty$ - guarantees
that $x$ be generic.\\

By Borel-Cantelli lemma ($\stackrel[k=1]{\infty}{\sum}\frac{1}{2^{k}}=1<\infty$),
$y\in\Lambda^{\mathbb{N}}$ $\mu$-almost-surely belongs - for $k$
large enough - to $\left[a_{k,r\left(y,k\right)}\right]$ for some $1\leq r\left(y,k\right)\leq p_{k}$. Denoting by $n\left(k,r\right)$
the $n$ for which $a'_{n}=a_{k,r}$, the fact that \\$a_{k,r\left(y,k\right)}=a'_{n\left(k,r\left(y,k\right)\right)}$
implies that $y$ agrees with $x$ on the indices \\$\left|a'_{n\left(k,r\left(y,k\right)\right)-1}\right|<i\leq\left|a'_{n\left(k,r\left(y,k\right)\right)-1}\right|+\left\lceil \frac{1}{3}\left|a'_{n\left(k,r\left(y,k\right)\right)}\right|\right\rceil $
and on \\
$\left|a'_{n\left(k,r\left(y,k\right)\right)-1}\right|+\left\lceil \frac{1}{3}\left|a'_{n\left(k,r\left(y,k\right)\right)}\right|\right\rceil <i\leq\left|a'_{n\left(k,r\left(y,k\right)\right)}\right|$
they differ by a shift of one symbol. Thus, recalling that $N_{k}>k\left|a_{k-1,p_{k-1}}\right|$,
the sequence in $k$ of $\sigma\times\sigma$-orbital measures of
$\left(x,y\right)$, taken for each $k$ until time \\
$\left\{ \left|a'_{n\left(k,r\left(y,k\right)\right)-1}\right|+\left\lceil \frac{1}{3}\left|a'_{n\left(k,r\left(y,k\right)\right)}\right|\right\rceil \right\} _{k\in\mathbb{N}}$
converges to the push-forward of $\mu$ through the map $z\mapsto\left(z,z\right)$
which let us denote by $\Delta_{\mu}$, while the similar sequence
taken for each $k$ until time $\left|a'_{n\left(k,r\left(y,k\right)\right)}\right|$
converges to the sum of $\frac{1}{3}\Delta_{\mu}$ and $\frac{2}{3}$
times the push-forward of $\mu$ through $z\mapsto\left(z,\sigma\left(z\right)\right)$.
So the point $\left(x,y\right)$ is not generic relative to any $\sigma\times\sigma$-invariant
measure. $\blacksquare$

Einstein Institute of Mathematics, Edmond J. Safra campus, The Hebrew
University of Jerusalem, Israel.\\
matan.tal@mail.huji.ac.il
\end{document}